\documentclass[11pt, oneside]{article}
\usepackage{amsmath,amsfonts,amssymb,amscd,theorem}
\date{}

\newtheorem{thm}{Theorem}
\newtheorem{prop}{Proposition}
\newtheorem{lem}{Lemma}

\begin{document}

\title{Meromorphic Differentials with Twisted Coefficients on Compact Riemann Surfaces}

\author{Yi-Hu Yang \footnote{Department of
Mathematics, Tongji University, Shanghai.} \thanks{The author
supported partially by NSF of China (No.10471105) and "Shuguang
Project" of Committee of Education of Shanghai (04SG21)}}

\maketitle
\begin{abstract} This note is to concern a generalization to the
case of twisted coefficients of the classical theory of Abelian
differentials on a compact Riemann surface. We apply the
Dirichlet's principle to a modified energy functional to show the
existence of differentials with twisted coefficients of the second
and third kinds under a suitable assumption on residues.
\end{abstract}

\section{Main results and discussion}
Let $\overline X$ be a compact Riemann surface. Classically, one
knows that a meromorphic (Abelian) differential can be expressed
as a sum of three kinds of differentials, one of which is
holomorphic, the second one {\it differentials of the second
kind}, i.e. all its poles having residues $0$, and the last one
{\it differentials of the third kind}, i.e. its poles being
$\log$-pole. A classical problem is, fixing some points in
$\overline X$, how to construct such a differential with poles at
the fixed points, {\it provided that the sum of residues be zero}.
This was completely solved, e.g. by using the Dirichlet's
principle on certain modified energy functional (cf. \cite{sie}).
Briefly, the results are as follows: For arbitrarily given point
$p$ of $\overline X$ with a local coordinate $z$ around $p$ and
arbitrary integer $k\ge 1$, one can find a differential $\phi$ of
the second kind such that $p$ is the only pole of $\phi$ and
$\phi$ has the following asymptotic behavior near $p$
\begin{equation}{\label{second kind differential}}
z^{-k-1}dz;
\end{equation}
for arbitrarily given two points $p_1, p_2$ of $\overline X$,
there exists a differential $\phi$ of the third kind such that
$p_1, p_2$ are the only $\log$-poles of $\phi$ and the residues of
$\phi$ are $1, -1$ at $p_1, p_2$ respectively; the general case
can be obtained by combining the above two. As mentioned above,
the method is the Dirichlet's principle; by using the Dirichlet's
principle on a certain modified energy functional, one can get a
harmonic function $u$ with prescribed asymptotic behaviors at the
given points and then $\partial u$ is the required Abelian
differential; {\bf the key is the requirement that the sum of
residues be zero.}

In this note, we want to generalize this classical theory to the
twisted case. Let $\rho: \pi_1(\overline X)\to Gl(n, \mathbb{C})$
be a linear representation of $\pi_1(\overline X)$, $L_\rho$ the
corresponding flat vector bundle, $D$ the canonical flat
connection on $L_\rho$. A Hermitian metric $h$ on $L_\rho$ can be
canonically explained as a $\rho$-equivariant map from the
universal covering of $\overline X$ into $Gl(n, \mathbb{C})/U(n)$
(equivalently, the set of all positive definite Hermitian
symmetric matrices, denoted by $\mathcal{P}_n$), still denoted by
$h$. Then, the differential $(dh)h^{-1}$ is a one-form valued in
${{End}}(L_\rho)$. The condition that the differential $(\partial
h)h^{-1}$ is holomorphic is then read as
\[
\overline\partial((\partial h)h^{-1})=0;
\]
equivalently, the map (metric) $h$ is harmonic (if $\overline X$
is higher dimensional, $h$ is pluri-harmonic). We consider
${{End}}(L_\rho)$ as our twisted coefficient. Then, our purpose of
this note is to find meromorphic one-forms with value in
${{End}}(L_\rho)$, which have prescribed singularities, similar to
classical Abelian differentials.

In order to find such differentials, we assume that the
representation $\rho: \pi_1(\overline X)\to GL(n, \mathbb{C})$ in
question is semi-simple (for the precise definition, see \S 3). We
attempt to find certain special $\rho$-eqivariant harmonic map
(harmonic metric) on $L_\rho$ with (possible) singularities;
equivalently, this means we apply the Dirichlet's Principle to
certain modified energy functional to get some special critical
points $h$ so that $(\partial h)h^{-1}$ are the desired ones. We
develop the variational technique of Siegel \cite{sie} so that it
is appropriate for the present nonlinear setup; this is one of
main points of this note.

In the following, we briefly describe our main results and their
proofs. Let me first show what our singularities look like. As in
the classical theory, we consider two kinds of singularities. We
first consider the second kind. Fix arbitrarily a point
$p\in\overline X$ and restrict ourself to a disk $\Delta$ with
center at $p$. Let $z=x+\sqrt{-1}y$ be an Euclidean complex
coordinate with $z(p)=0$. Restrict the flat bundle $L_\rho$ to
$\Delta$ and fix a suitable flat basis of $L_\rho$ on $\Delta$.
Then, under the fixed basis, the asymptotic behavior
{\footnote{for the precise definition of the asymptotic behavior
of a metric at a puncture, see \S4.}} at the point $p$ of the
desired harmonic metrics is of the following form
\begin{equation}\label{asymp-behavior}
{\exp}\left(
\begin{array}{cccc}
\sum a_{k_1}u_{k_1}&~&0\\
~&\ddots&~\\
0&~&\sum a_{k_n}u_{k_n}
\end{array}
\right)
\end{equation}
where $k_1, k_2, \cdots, k_n\in \mathbb{N}$, $a_{k_1}, \cdots,
a_{k_n}\in \mathbb{R}$, and $u_k=2{\text{Re}}(z^{-k})$. It is easy
to see that if a harmonic metric $K$ has the above asymptotic
behavior, the corresponding differential $(\partial K)K^{-1}$ then
has the following asymptotic behavior
\begin{equation}\label{higgs-asymp-behavior}
\left(
\begin{array}{cccc}
-\sum k_1a_{k_1}{z^{-k_1}}&~&0\\
~&\ddots&~\\
0&~&-\sum k_na_{k_n}{z^{-k_n}}
\end{array}
\right)\cdot{\frac{dz}{z}},
\end{equation}
Our first result is then the following

\begin{thm}  Let $\rho: \pi_1(\overline X)\to Gl(n,
\mathbb{C})$ be a semi-simple representation (for the precise
definition, see \S 3), $p_1, p_2, \cdots, p_s$ arbitrarily given
points of $\overline X$; by $X$ denote $\overline X\setminus\{p_1,
p_2, \cdots, p_s\}$. Let $L_\rho$ be the corresponding flat bundle
restricted to $X$. Then, for arbitrarily given asymptotic
behaviors of the form \label{asymp-behavior}(2) at the punctures
$p_1, p_2, \cdots, p_s$, there exists a unique harmonic metric $K$
on $L_\rho$ with the corresponding asymptotic behaviors; hence the
differential $(\partial K)K^{-1}$ is a holomorphic one-form with
twisted coefficient which has asymptotic behavior of the form (3).
\end{thm}

We now consider the singularities of the third kind. Fix
arbitrarily points $p_1, p_2, \cdots, p_s\in\overline X$, take a
smooth curve $\gamma$ connecting them, say, the starting point
$p_1$, the end point $p_s$; take a small enough tube neighborhood
$\Gamma$ of $\gamma$ so that they are simply-connected. Now,
assume that under a fixed flat basis of $L_\rho$ on $\Gamma$, the
desired differentials at each $p_i$ have prescribed singularity of
the following form
\begin{equation}\label{log-asymp-behavior}
\left(
\begin{array}{cccc}
a_1^i&~&0\\
~&\ddots&~\\
0&~&a_n^i
\end{array}
\right)\cdot{\frac{dz^i}{z^i}},
\end{equation}
where $z^i$ is a local complex coordinate at $p_i$ and $a_1^i,
\cdots, a_n^i\in\mathbb{R}$. Then, we have the following

\begin{thm} Let $\rho$ as before, $p_1, \cdots, p_s$ arbitrarily fixed points on
$\overline X$; by $X$ denote $\overline X\setminus\{p_1, \cdots,
p_s\}$. Let $L_\rho$ be the corresponding flat bundle restricted
to $X$. Then, for arbitrarily given asymptotic behaviors of the
form \label{asymp-behavior1}(4) at the points $p_i$ such that the
$a^i_j$ are rational numbers (actually, we can assume the ratios
of $a_j^i$ and $a_j^k$ are rational; see \S3) and
$\sum_{i=1}^sa_j^i=0, j=1, \cdots, n$, there exists a unique
harmonic metric $K$ on $L_\rho$, the differential of which
$(\partial K)K^{-1}$ is a holomorphic one-form with twisted
coefficient and has asymptotic behavior of the form (4) at each
point $p_i$.
\end{thm}

We now outline the proof of the theorems. The proofs of Theorem 1
and 2 are completely similar, so we outline only that of Theorem
1. It is clear that the harmonic metrics with prescibed
singularities are always of infinite energy. So, the variational
technique for the usual energy functional does not work anymore.
In order to overcome this difficulty, we use a modified energy
functional, which is roughly defined as follows. Let $\mathcal{K}$
be the set of continuous and piece-wise differentiable metrics on
$L_\rho$ which have the asymptotic behaviors mentioned above near
each puncture $p_i$. For a metric $K\in\mathcal{K}$, we define its
modified energy as
\begin{equation}\label{modified-energy}
\hat{E}(K)=\int_{X\setminus\bigcup_{i=1}^s\Delta_{i}^*}\mid
(dK)K^{-1}\mid^2+\sum_{i=1}^s\int_{\Delta_{i}^*}\mid (dK)K^{-1}-
(dK_0)K_0^{-1}\mid^2,
\end{equation}
where $K_0$ is a suitably constructed metric with the asymptotic
behaviors above at each puncture $p_i$, $\Delta^*_i$ is a small
punctured disk around $p_i$. We remark that, in the definition of
the modified energy, we use the difference of the derivatives of
two maps, which applies to the only case when the target manifolds
are homogeneous. Then, we will prove that one can minimize the
modified energy functional $\hat{E}(K)$ in $\mathcal{K}$ and the
minimizer is a (smooth) harmonic metric with prescribed asymptotic
behaviors at the punctures.

In order to minimize the modified energy functional $\hat{E}$,
technically, we first need to construct a suitable initial metric
$K_0$, which is harmonic near each puncture $p_i$ and has not only
prescribed asymptotic behaviors at the punctures but also
vanishing radial derivatives on a certain circle around each
puncture. We would like to point out that both the harmonic
property of $K_0$ around the punctures and vanishing radial
derivatives of $K_0$ on a certain circle around each puncture are
very key for our proof. After this construction, we choose a
minimizing sequence of $\hat{E}$ in $\mathcal{K}$. Generally, such
a minimizing sequence does not necessarily converge. In order to
make such a sequence to converge, we have to modify it. To this
end, we first use harmonic metrics to replace continuously each
metric of the minimizing sequence on $X\setminus\cup_i\Delta_i^*$.
It is clear that the new sequence is still a minimizing one of
$\hat{E}$; furthermore, using the semi-simplicity of $\rho$, we
can show that the new sequence (if necessary, going to a
subsequence) is actually uniformly convergent on any compact
subset of $X\setminus\cup_i\Delta_i^*$. We continue to modify the
new minimizing sequence on the remaining part. For this, we need
to solve a boundary value problem for harmonic metrics with
prescribed asymptotic behavior at the puncture on a punctured disk
(Proposition 4). After solving such a boundary problem, we then
use such a solution to replace continuously each metric of the new
minimizing sequence on a greater disk than $\Delta^*_i$ around
each puncture; we can show that the sequence obtained is still a
minimizing one (Proposition 5). Using the previous convergence, we
can finally show that the sequence obtained is uniformly
convergent on $X$ and the limit lies in $\mathcal{K}$.

From the above description, it is easy to see that our proof for
convergence of minimizing sequences is slightly different from
that of Siegel; we use a two-step modification of minimizing
sequences and the semi-simplicity of the representation $\rho$. In
fact, although the argument for convergence of Siegel can be
explained as the case of one-dimensional trivial representations,
it, due to the nonlinearity of maps, does not however apply to the
present setting.

The idea of modifying energy functional was also used by Ding in
\cite{di} to deal with the problems of harmonic maps with infinite
energy. Due to generality of the target manifolds he considered,
he can not use the difference of the derivatives of two maps;
instead, he used the integration by parts on bounded domains of
the domain manifolds and then an approximation process.

Naturally, one should ask if there exists a holomorphic one-form
with twisted coefficients but without singularity. In the case of
complex coefficient, this is a well-known result; the dimension of
the set of such differentials is the genus of $\overline X$. In
the case of twisted coefficient, this is actually a consequence of
Donaldson's result \cite{do}; in the case of higher dimension,
this is also true by means of Corlette's result \cite{co} and
Siu's Bochner technique. In a future paper, we will generalize the
results and the method of the present paper to the higher
dimension case.

Some of the points of the present work were realized when I was
visiting The University of Hong Kong during June-July of 2005,
where Professor Yum-Tong Siu explained to me the idea of Siegel
and its physical background; I thank him very much for kindly
letting me share his idea and for later several talks which are
important for my understanding of the problems here. I also thank
Professor Ngming Mok for inviting me to visit HKU. Thanks also
goes to Yuxin Ge for some helpful discussion about the proof of
Proposition 4. Part of the paper was done when I stayed in
Max-Planck Institute for Mathematics, Leipzig as a guest; I would
like to thank Professor J\"urgen Jost for his support and the
institute for its hospitality and a good working environment. Part
of the results was talked in The Second Sino-German Conference on
Complex Geometry, Shanghai, September 2006.

\section{The energy functional and the equation}
In this section, for convenience, we fix some notations and state
some more or less standard facts (cf. e.g. \cite{si1, si2}). Let
$\mathcal{P}_n$ be the set of all positive definite hermitian
symmetric matrices of order $n$. $GL(n, \mathbb{C})$ acts
transitively on $\mathcal{P}_n$ by
\[
g\circ A={g}A^t\bar{g}, ~A\in\mathcal{P}_n, g\in GL(n,
\mathbb{C}).
\]
Obviously, the action has the isotropic subgroup $U(n)$ at the
identity $I_n$. Thus $\mathcal{P}_n$ can be identified with the
coset space $GL(n, \mathbb{C})/U(n)$, and can be uniquely endowed
an invariant metric{\footnote{In terms of matrices, such an
invariant metric can be defined as follows. At the identity $I_n$,
the tangent elements just are hermitian matrices; let $A, B$ be
such matrices, then the Riemannian inner product $<A,
B>_{\mathcal{P}_n}$ is defined by $tr(AB)$. In general, let $H\in
\mathcal{P}_n$, $A, B$ two tangent elements at $H$, then the
Riemannian inner product $<A, B>_{\mathcal{P}_n}$ is defined by
$tr(AH^{-1}BH^{-1})$.} up to some constants. In particular, under
such a metric, the geodesics through the identity $I_n$ are of the
form $\exp(th)$, $t\in \mathbb{R}$, $h$ being a hermitian matrix.

Let $X$ be a complex manifold, $\mathbb{V}\to X$ a flat vector
bundle, $K$ a hermitian metric on $\mathbb{V}$. For $x\in X$, the
metric $K_x$ on the fiber $\mathbb{V}_x$, after fixing a basis,
can be considered as an element $H_x\in \mathcal{P}_n$, and hence
a point in the coset space $GL(n, \mathbb{C})/U(n)$. Thus, after
fixing a flat basis of $\mathbb{V}$, the metric $K$ can be
considered as an equivariant map from the universal covering of
$X$ into $GL(n, \mathbb{C})/U(n)$ or $\mathcal{P}_n$.

From now on, we always fix a flat basis $\{v_i\}$ of $\mathbb{V}$,
unless stated otherwise. Decompose the flat connection $D=d'+d''$
into the parts of type $(1,0)$ and $(0,1)$. Define the
differential operators $\delta'$ and $\delta''$ by setting
\begin{eqnarray*}
&&\partial<u, v>_K=<\delta' u, v>_K+<u, d''u>_K,\\
&&\overline\partial<u, v>_K=<\delta'' u, v>_K+<u, d'u>_K,
\end{eqnarray*}
namely, both $\delta'+d''$ and $d'+\delta''$, as connection on
$\mathbb{V}$, preserve the metric. Clearly $d''d'+d'd''=0$ implies
$\delta'\delta''+\delta''\delta'=0$. Set
\begin{eqnarray*}
\theta_K&=&(d'-\delta')/2, \overline\theta_K=(d''-\delta'')/2,\\
\partial_K&=&(d'+\delta')/2, \overline\partial_K=(d''+\delta'')/2.
\end{eqnarray*}
It is easy to see that $<\theta_K u, v>=<u, \overline\theta_K v>$
and $\partial_K+\overline\partial_K$ preserves the metric.
$\theta_K$ (resp. $\overline\theta_K$) is a one-form of type $(1,
0)$ (resp. $(0, 1)$) valued in $End(\mathbb{V})$
{\footnote{Actually, the construction of the operators
$\partial_K, \overline\partial_K, \theta_K, \overline\theta_K$
comes essentially from the Cartan decomposition of the flat
connection $D$ with respect to $\mathfrak{gl}(n,
\mathbb{C})=\mathfrak{u}(n)+\mathfrak{p}_n$, where
$\mathfrak{p}_n$ is the set of Hermitian matrices of order $n$.
This, together with the fact that $\mathcal{P}_n$ is homogeneous,
implies that these operators are invariant under certain sense, as
will be used in various computations of the note very often. The
point can be more explicitly seen if we consider $K$ as an
equivariant map into $\mathcal{P}_n$, so that the connection
$\partial_K+\overline\partial_K$ is the pull-back of the standard
invariant connection of $\mathcal{P}_n$}}. On the other hand, one
can explicitly write down $\theta_K$ in terms of the basis
$\{v_i\}$ as follows. Setting $H=(H_{i\bar j})=(<v_i, v_j>_K)$,
one has then
$$
\partial H_{i\bar j}=\partial
<v_i, v_j>_K=<\delta' v_i, v_j>_K=-2<\theta_K v_i, v_j>_K.
$$
Writing $\theta_K=\theta^k_{i\alpha}v_k\otimes v^i\otimes
dz^\alpha$ ($\{v^i\}$ is the dual basis of $\{v_i\}$,
$\{z^\alpha\}$ is a local coordinate of $X$), one then has
\[
\partial H_{i\bar j}=-2\theta^k_{i\alpha}H_{k\bar j}dz^\alpha,
\]
namely,  $\partial H_{i\bar k}H^{\bar k
j}=-2\theta^j_{i\alpha}dz^\alpha$ (or invariantly, $\partial
HH^{-1}=-2\theta_K$), where $(H^{\bar k j})$ is the inverse of
$(H_{i\bar j})$. Thus, $\theta_K$ (resp. $\overline\theta_K$) can
be identified with the differential of type $(1, 0)$ (resp. $(0,
1)$) of the map into $\mathcal{P}_n$ corresponding to the metric
$K$, up to some constant.
\\
\\
%\vskip .2cm \noindent
{\it Remark.} $\delta'+d''=D-2\theta_K$ can
also be regarded as a (hermitian) connection on the flat
($d''$-holomorphic) bundle $\mathbb{V}$ with respect to $K$, so
that $-2\theta_K$ is just the corresponding connection form under
the fixed flat basis, as will be used in the following. The
similar explanation works for $d'+\delta''$ and
$\overline\theta_K$.
\\

%\vskip .2cm
Usually, one needs to choose some "nice" metrics of
$\mathbb{V}$, which furthermore satisfy some differential
equations. To this end, we from now on assume that $X$ is a
K\"ahler manifold with a K\"ahler metric $\omega$, and denote by
$\Lambda$ the adjoint of the operation of wedging with $\omega$ on
exterior forms of $X$. Set
$$
D''_K=\overline\partial_K+\theta_K,
~~D'_K=\partial_K+\overline\theta_K
$$
and $G_K=(D''_K)^2$; call $D_K''$ the {\it Higgs operator} and
$G_K$ the {\it pseudo-curvature}. Call a metric $K$ on the flat
bundle $\mathbb{V}$ {\it harmonic} if it satisfies
\begin{equation}\label{harmoni-equation}
\Lambda G_K=0.
\end{equation}
Equivalently, this can be written as
$$
\Lambda(d'\delta''+\delta''d'-\delta'd''-d''\delta')=0.
$$
In the following, we will show that the metric $K$ being harmonic
is equivalent to the corresponding map being a harmonic map, and
hence {\it the equation $\Lambda G_K=0$ is a variational one}.

First, let us see how $\theta_K$ and $\overline\theta_K$ change
when the metric $K$ changes. Let $K_1, K_0$ be two metrics on
$\mathbb{V}$. One can then define an endomorphism $h$ of
$\mathbb{V}$ by setting $<u, v>_{K_1}=<hu, v>_{K_0}$. It is clear
that $h$ is self-adjoint positive with respect to $K_0$. Under the
fixed basis $\{v_i\}$, write $h=(h_i^j)$, i.e. $hv_i=h_i^jv_j$;
also write $K_s$ as the matrix $(H_{si\bar j})$, the inverse of
which is denoted by $H_s^{~i\bar j}$, $s=0, 1$. Then $H_{1i\bar
j}=h_i^kH_{0k\bar j}$ and $H_1^{~\bar ji}=H_0^{\bar
jk}(h^{-1})_k^i$, here $h^{-1}$ is the inverse of $h$. Thus
\begin{eqnarray*}
-2\theta_{K_1}&=&\partial H_1H_1^{-1}\\
&=&(\partial hH_0+h\partial H_0)H_0^{-1}h^{-1}\\
&=&\partial hh^{-1}+h\partial H_0H^{-1}_0h^{-1}\\
&=&\partial hh^{-1}-2h\theta_{K_0}h^{-1} \\
&=&(\partial h-2h\theta_{K_0}+2\theta_{K_0}h-2\theta_{K_0}h)h^{-1} \\
&=&\delta'hh^{-1}-2\theta_{K_0},
\end{eqnarray*}
in the last equality, $h$ is considered as a section of
$End(\mathbb{V})$, and $\delta'h$ is the covariant derivative of
$h$ (referring to the above remark). Similarly,
\begin{equation}\label{theta-transformation}
\overline\theta_{K_1}=-(1/2)\delta''hh^{-1}+\overline\theta_{K_0}.
\end{equation}
The above computation is very similar to that in the
Hermitian-Yang-Mills theory (cf. e.g. \cite{siu}).

For later purpose, let's here make some simple remarks about both
inner products on the bundle $End{\mathbb{V}}$ and the tangent
bundle of $\mathcal{P}_n$. In the following arguments, we ignore
integrability of integrals. For $End{\mathbb{V}}$, we always use
the trace inner product $<\mathcal{A},
\mathcal{B}>=tr(\mathcal{A}^t\overline{\mathcal{B}})$ under the
fixed flat basis; it is easy to see that when endomorphisms are
from tangent elements of $\mathcal{P}_n$, the trace inner product
coincides with the invariant inner product on $\mathcal{P}_n$. Fix
a metric $K$ on $\mathbb{V}$ (and hence $End(\mathbb{V})$) and a
point $x\in X$. The metric $K_x$ on the fiber $\mathbb{V}_x$
corresponds to the matrix $H_x\in\mathcal{P}_n$. Let $A, B$ be two
elements in the tangent space $T_{H_x}\mathcal{P}_n$. Then
$AH^{-1}_x, BH^{-1}_x$ can be considered as two self-adjoint
homomorphisms of $\mathbb{V}_x$ with respect to $K_x$. Using an
orthogonal basis of $\mathbb{V}_x$ with respect to $K_x$, it is
easy to show that
\[
<A, B>_{\mathcal{P}_n}=tr(AH^{-1}_xBH^{-1}_x)=<AH^{-1}_x,
BH^{-1}_x>_{K_x}.
\]
This together with the previous argument concerning $\theta_K$ and
$\overline\theta_K$ implies
\begin{eqnarray*}
&&\int_X(<\theta_KH, \theta_KH>_{\mathcal{P}_n,
\omega}+<\overline\theta_KH, \overline\theta_KH>_{\mathcal{P}_n,
\omega})\\
&=&\int_X(<\theta_K, \theta_K>_{K, \omega}+<\overline\theta_K,
\overline\theta_K>_{K, \omega}).
\end{eqnarray*}
Based on these remarks, afterwards, we often omit the subscripts
of the inner products, since it is clear from the context.

Furthermore, both integrals above are independent of choice of a
basis, though defined by choosing a basis; and the first integral,
up to a constant, is just the energy of the map $H$ corresponding
to the metric $K$; for simplicity, we call it the energy of the
metric $K$, denoted by $E(K)$ , i.e.,
\begin{equation}\label{metric-energy}
E(K)=\int_X(<\theta_K, \theta_K>+<\overline\theta_K,
\overline\theta_K>).
\end{equation}
Thus, taking the first variation for either of both integrals with
respect to $K$, we will get the equation of harmonic maps for $K$.
As usual, we do this for the first integral and, due to
$\mathcal{P}_n$'s homogeneity, we can easily pass to an integral
related to the bundle $\mathbb{V}$. Finally the Euler-Lagrange
equation is just the equation of harmonic metrics $\Lambda G_K=0$.

Let $h$ be a self-adjoint (not necessarily positive) endomorphism
of $\mathbb{V}$ with respect to $K$. Then $\exp(th)$ is
self-adjoint positive, $t\in \mathbb{R}$. Set $H_t=\exp(th)H$,
denote the corresponding metric by $K_t$, here $H=(H_{i\bar
j})=<v_i, v_j>_K=H_0$. From the previous computation, one has
\begin{eqnarray*}
-2\theta_{K_t}&=&\delta'(\exp(th))\exp(-th)-2\theta_{K}\\
&=&t\delta'h-2\theta_K+o(t).
\end{eqnarray*}
Similarly,
$-2\overline\theta_{K_t}=t\delta''h-2\overline\theta_K+o(t)$.
Thus,
\begin{eqnarray*}
{\frac{d}{dt}}E(K_t)_{|_{t=0}}&=&{\frac{d}{dt}}\int_X(<\theta_{K_t}H_t,
\theta_{K_t}H_t>_{\mathcal{P}_n,
\omega}+<\overline\theta_{K_t}H_t,
\overline\theta_{K_t}H_t>_{\mathcal{P}_n,
\omega})_{|_{t=0}}\\
&=&-\int_X(<\theta_{K}H, (\delta'h)H>_{\mathcal{P}_n,
\omega}+<\overline\theta_{K}H,
(\delta''h)H>_{\mathcal{P}_n, \omega})\\
&=&-\int_X(<\theta_{K}, \delta'h>_{K,
\omega}+<\overline\theta_{K},
\delta''h>_{K, \omega})\\
&=&-\int_X(<(\delta')^*\theta_{K}, h>_{K,
\omega}+<(\delta'')^*\overline\theta_{K}, h>_{K, \omega}).
\end{eqnarray*}
Since $\delta'+d''$ (resp. $d'+\delta''$) is a hermitian
connection on $\mathbb{V}$ with respect to $d''$ (resp.
$\delta''$) and $K$, so one has the related K\"ahler identity
$(\delta')^*=\sqrt{-1}[\Lambda, d'']$ (resp.
$(\delta'')^*=-\sqrt{-1}[\Lambda, d']$). Thus,
\begin{eqnarray*}
{\frac{d}{dt}}E(K_t)|_{t=0}&=&-\int_X\sqrt{-1}(<\Lambda
d''(\theta_{K}), h>_{K, \omega}-<\Lambda d'(\overline\theta_{K}),
h>_{K, \omega})\\
&=&-\int_X\sqrt{-1}<\Lambda
(d''(\theta_{K})-d'(\overline\theta_{K})), h>_{K, \omega}.
\end{eqnarray*}
So, the E-L equation is
\begin{equation}\label{harmonic-equation-2}
\sqrt{-1}\Lambda (d''(\theta_{K})-d'(\overline\theta_{K}))=0,
\end{equation}
which is just equivalent to the equation $\Lambda G_K=0$ for $K$
being harmonic. Summing up all the above argument, we have
\begin{prop}
Let $(X, \omega)$ be a K\"ahler manifold, $\mathbb{V}\to X$ a flat
vector bundle. Giving a metric $K$ on $\mathbb{V}$ and using the
above notations, one has the energy functional
\[
E(K)=\int_X\mid\theta_K+\overline\theta_K\mid^2,
\]
the E-L equation of which is
$$
\sqrt{-1}\Lambda (d''(\theta_{K})-d'(\overline\theta_{K}))=0,
$$
i.e. $\Lambda G_K=0$.
\end{prop}

\section{The construction of initial metrics}
Let $X$ be a differentiable manifold. In general, we call a linear
representations $\rho: \pi_1(X)\to Gl(n, \mathbb{C})$ semi-simple
if the Zariski closure in $Gl(n, \mathbb{C})$ of the image of
$\rho$, as an algebraic subgroup, is semi-simple. Here, for
convenience of the later application, we state a more geometric
definition of semi-simplicity, which is motivated by Donaldson
\cite{do}; for this, we need to use a little bit knowledge about
the boundary theory of symmetric spaces (cf. e.g. \cite{bj}).

Let $\rho: \pi_1(X)\to Gl(n, \mathbb{C})$ be a linear
representation. Call $\rho$ {\it semi-simple} if for any boundary
component $\Sigma$ of $\mathcal{P}_n$, there exists an image
element $\rho(\gamma), \gamma\in\pi_1(X)$ satisfying
$\Sigma\cap\rho(\gamma)(\Sigma)=\emptyset$.

From now on, we assume that $X$ is an open Riemann surface, i.e.
$X=\overline X\setminus\{p_1, p, \cdots, p_s\}$ for a compact
Riemann surface $\overline X$, as mentioned in the \S1. Let $\rho:
\pi_1(\overline X)\to Gl(n, \mathbb{C})$ be a linear semi-simple
representation. Let $L_\rho$ be the corresponding flat bundle over
$\overline X$; also by $L_\rho$ denote the restriction to $X$.

We first construct an initial metric on $L_\rho$ needed by the
proof of Theorem 1. Let $\Delta^*_i$ (resp. $\Delta^*_{i/2}$) be
the punctured disk with radius $1$ (resp. $\frac 1 2$) around
$p_i$ and $t_i (=r_ie^{\sqrt{-1}\theta})$ a complex (polar)
coordinate on $\Delta^*_i$ with $t_i(p_i)=0$. Fix a flat basis of
$L_\rho$ which is clearly single-value on each $\Delta^*_i$. Under
this flat basis, we construct a metric of $L_\rho$ on $\Delta^*_i$
as follows:
\begin{equation}\label{asymp-behavior-1}
H_i={\exp}\left(
\begin{array}{cccc}
\sum a_{k_1}u_{k_1}&~&0\\
~&\ddots&~\\
0&~&\sum a_{k_n}u_{k_n}
\end{array}
\right),
\end{equation}
where  $k_1, k_2, \cdots, k_n\in\mathbb{N}$, $a_{k_1}, a_{k_2},
\cdots, a_{k_n}\in\mathbb{R}$, and
$$
u_k={\text{Re}}({\frac{1}{t^k_i}}+4^kt_i^k).
$$
It is clear that the $H_i$ is a harmonic metric on $\Delta^*_i$
and
\[
{\frac{\partial H_i}{\partial r_i}}H_i^{-1}=0,~~ {\text{on}}
~r_i={\frac 1 2}.
\]
Extending the $H_i$'s to $X$, we obtain the required initial
metric of $L_\rho$, denoted by $K_0$.

We now turn to the construction of the initial metrics needed in
the proof of Theorem 2. As in \S1, connect $p_1, \cdots, p_s$ in a
smooth curve $\gamma$ and take two small enough neighborhoods
$\Gamma\subset\Gamma'$ of $\gamma$ so that they are
simply-connected. Take a compatible (with the complex structure of
$\overline X$) complex coordinate $z$ on $\Gamma'$ so that
$\Gamma$ and $\Gamma'$ are two disks; without loss of generality,
assume $z(p_1)=0$. Put $\Gamma'$ on the complex plane and take the
reflections of $p_2, \cdots, p_s$, denote by $p_2', \cdots, p_s'$
respectively. (If necessary, we can shrink $\Gamma'$ so that the
reflection points do not lie in $\Gamma'$.) Denote the coordinates
of $p_2, p_2', \cdots, p_s, p_s'$ by $\xi_2, \xi_2', \cdots,
\xi_s, \xi_s'$ respectively. Construct a meromorphic function on
$\Gamma'$ as follows
\[
g(z)={\frac{z^{l_2+\cdots
+l_s}}{\prod_{i=2}^s(z-\xi_i)^{l_i}\cdot\prod_{i=2}^s(z-\xi_i')^{l_i}}},
\]
where $l_2, \cdots, l_s\in\mathbb{N}$; take the real part of the
multiple-value function $\log g(z)$, which is single-value,
denoted by $u_{l_2\cdots l_s}$. A simple argument shows that the
radial derivatives of $u_{l_2\cdots l_s}$ along $\partial\Gamma$
vanish. Then, under a fixed flat basis of $L_\rho$ on $\Gamma'$,
our initial metric on $L_\rho$ over $\Gamma'$ is taken as the
following
\begin{equation}\label{asymp-behavior1}
H={\exp}\left(
\begin{array}{cccc}
a_1u_{l_2^1\cdots l_s^1}&~&0\\
~&\ddots&~\\
0&~&a_nu_{l_2^n\cdots l_s^n}
\end{array}
\right),
\end{equation}
where $a_1, \cdots, a_n\in \mathbb{R}$. It is clear that $H$ is a
harmonic metric on $\Gamma'\setminus{p_1, p_2, \cdots, p_s}$ and
\[
{\frac{\partial H}{\partial r}}H^{-1}=0,~~ {\text{on}}~~
\partial\Gamma,
\]
where $r$ is the radial coordinate of $z$. Now, we can extend $H$
to $X$ to get a desired initial metric on $L_\rho$ with prescribed
singularities, also denoted by $K_0$. It is also clear that If a
harmonic metric $K$ has the above asymptotic behavior at the
points $p_i$, the corresponding differential $(\partial K)K^{-1}$,
under the fixed flat basis, has then asymptotic behavior of the
following form at the points $p_i$
\begin{equation}\label{higgs-asymp-behavior1}
\left(\begin{array}{cccc}
b_1&~&0\\
~&\ddots&~\\
0&~&b_n
\end{array}
\right)\cdot{\frac{dz}{z}},
\end{equation}
where $b_1, \cdots, b_n\in\mathbb{R}$.

\section{The modified energy functional and minimizing sequences}
In this section, we develop a variational technique to show the
existence of a harmonic metric on $L_{\rho}$ with the prescribed
asymptotic behaviors at the punctures $p_i$, as described in \S3;
since the proofs of Theorem 1 and 2 are the same, our discussion
here is restricted to Theorem 1. In order to use the method of
minimizing sequences to get such a harmonic metric, we need to
modify the usual energy functional $E(K)$. Due to conformal
invariance, without loss of generality, we can take a special
Riemannian metric $\omega$ on $X$ which is Euclidean on each
$\Delta_i^*$, i.e., $\omega_{|\Delta^*_i}=\sqrt{-1}dt_i\wedge
d\overline{t}_i$.
\\

Call a (continuous and piece-wise differentiable) metric $K$
($=hK_0$ under a fixed flat basis) on $L_\rho$ (or
${L_\rho}_{|\cup_{i=1}^s\Delta^*_i}$) {\it admissible} relative to
$K_0$ if
\\

1)~~The integral
$\sum_{i=1}^s\int_{\Delta_{i/2}^*}\big(\mid\theta_{K}-
\theta_{K_0}\mid^2+\mid\overline\theta_{K}-\overline\theta_{K_0}\mid^2\big)
$ exists; and

2)~~$K$ and $K_0$ are mutually bounded (namely, if writing
$K=hK_0$, the eigenvalues of $h$ are uniformly far away from both
$0$ and $\infty$).
\\

We remark that when one considers the metrics on $L_\rho$ as
equivariant maps from the universal covering $\tilde X$ of $X$
into $\mathcal{P}_n$, the condition 2) is equivalent to
\\

2)'~~$K$ and $K_0$ have uniform bounded distance near the
punctures under the invariant metric of $\mathcal{P}_n$.
\\

The equivalence of 2) and 2)' will be used very often in the
following discussion. If two metrics $K_1, K_2$ satisfy the
property 2) or 2)', we say they have the same asymptotic behavior
at the punctures.
\\

Denote the set of admissible metrics $K$ on $L_{\rho}$ by
$\mathcal{K}$. We then define the {\it modified energy functional}
as
\begin{equation}\label{4.1}
\hat{E}(K)=\int_{X\setminus\bigcup_{i=1}^s\Delta_{i/2}^*}\mid\theta_K+
\overline\theta_K\mid^2+\sum_{i=1}^s\int_{\Delta_{i/
2}^*}\big(\mid\theta_{K}-
\theta_{K_0}\mid^2+\mid\overline\theta_{K}-\overline\theta_{K_0}\mid^2\big),
\end{equation}
for admissible hermitian metrics $K\in\mathcal{K}$ on $L_{\rho}$.
For the modified energy functional, we have

\begin{prop}
Any admissible metric $K$ on $L_\rho$, if it minimizes the
modified energy functional $\hat{E}$, is a harmonic metric (and
hence smooth).
\end{prop}
{\it Proof.} One only needs to use the same computation as in \S 2
(taking the first variation of $\hat E$ to obtain the E-L
equation) and remarks that $K_0$ is harmonic on
$\cup_{i=1}^s\Delta^*_i$ and has vanishing normal derivative on
$\partial\Delta_{i/2}$. Mainly, one needs to consider variational
domains containing (part of) $\partial\Delta_{i/2}$, then the
condition of vanishing normal derivative applies; since if a
variational domain lies completely in the interior of
${X\setminus\bigcup_{i=1}^s\Delta_{i/2}^*}$ or
$\bigcup_{i=1}^s\Delta_{i/2}^*$, the minimizer $K$ of $\hat{E}$ is
naturally harmonic in such a domain using the harmonicity of $K_0$
on each $\Delta^*_i$ and the usual first variation computation.

Suppose that $D$ is a compact sub-domain containing (part of)
$\partial\Delta_{i/2}$. Consider a one-parameter variation $K^t$
of $K$ with $K^0=K$ and $K^t_{|X\setminus D}\equiv K$;
correspondingly, $K^t=h^tK$ with $h^0=$identity and
$h^t_{|X\setminus D}\equiv$identity; furthermore, we can even
assume that $h^t=e^{th}$ with $h$ being hermitian and
$h_{|X\setminus D}\equiv 0$. For sake of simplicity, by $X^c$
denote ${X\setminus\bigcup_{i=1}^s\Delta_{i/2}^*}$. Compute at
$t=0$
\begin{eqnarray*}
&&{\frac{d}{dt}}_{|t=0}\hat{E}(K_t)\\
&=&\int_{X^c}{\frac{d}{dt}}\mid\theta_{K^t}+
\overline\theta_{K^t}\mid^2+\sum_{i=1}^s\int_{\Delta_{i/2}^*}{\frac{d}{dt}}\big(\mid\theta_{K^t}-
\theta_{K_0}\mid^2+\mid\overline\theta_{K^t}-\overline\theta_{K_0}\mid^2\big)\\
&=&\int_{X^c\cap D}<\theta_{K}+ \overline\theta_{K}, \delta'
h+\delta''h>+\sum_{i=1}^s\int_{\Delta_{i/2}^*\cap D}<\theta_{K}-
\theta_{K_0}+\overline\theta_{K}-\overline\theta_{K_0}, \delta'h+\delta''h >\\
&=&\int_{D}<\theta_{K}+ \overline\theta_{K}, \delta'
h+\delta''h>-\sum_{i=1}^s\int_{\Delta_{i/2}^*\cap D}<\theta_{K_0}+\overline\theta_{K_0}, \delta'h+\delta''h >\\
&=&\int_{D}<\theta_{K}+ \overline\theta_{K}, \delta'
h+\delta''h>+\sum_{i=1}^s\int_{\Delta_{i/2}^*\cap
D}d<\theta_{K_0}+\overline{\theta}_{K_0},
h>+\\
&&+\sum_{i=1}^s\int_{\Delta_{i/2}^*\cap
D}<\Lambda(d''\theta_{K_0}-d'\overline{\theta}_{K_0}), h>\\
&=&\int_{D}<\theta_{K}+ \overline\theta_{K}, \delta' h+\delta''h>,
\end{eqnarray*}
in the last equality, we use the Stokes' formula, the harmonic
property of $K_0$ on each $\Delta^*_i$, vanishing normal
derivative of $K_0$ on $\partial\Delta_{i/2}$ with respect to
$\Delta^*_{i/2}$, and $h$ being vanishing on $\partial D$.

Finally, the minimizing property of $K$ for $\hat{E}$ implies that
${\frac{d}{dt}}_{|t=0}\hat{E}(K_t)$=0. On the other hand, since
$h$ can be chosen arbitrarily with $h_{|X\setminus D}\equiv 0$,
so, under the weak sense,
\[
\Lambda(d''\theta_K-d'\overline{\theta}_K)=0,
\]
namely, $K$ is weakly harmonic on $D$, and hence harmonic.
~~~~~~~~~~~$\Box$
\\

Since $\hat{E}(K)\ge 0$ for any $K\in \mathcal{K}$, so the
greatest lower bound of $\hat{E}$ on $\mathcal{K}$ is a
nonnegative number, denoted by $\mu$. Therefore, there exists a
sequence of admissible metrics $\{K_n\}_{n=1}^\infty$ with
$\lim_{n\to \infty}{\hat E}(K_n)=\mu$. Call such a sequence a {\it
minimizing sequence} of $\hat{E}$ in $\mathcal{K}$. In general, it
is not clear if such a minimizing sequence is convergent and (if
so) the limit is an admissible harmonic metric. We shall however
show that it is possible to prove that minimizing sequences
constructed in a special way are convergent and the corresponding
limits are admissible and harmonic, so that the greatest lower
bound $\mu$ is actually attained as the minimum value of $\hat E$
for an admissible metric which has the same behavior as $H_i$ at
each $p_i$. Namely, we will show
\\
\\
{\bf Main Theorem} {\it There exists an admissible metric $K\in
\mathcal{K}$ with ${\hat E}(K)=\mu$, and hence $K$ is harmonic.}
\\

In order to prove the main theorem, we here recall some estimates
for harmonic maps into non-positive curved manifolds, which are
presently standard and also apply very well to the equivariant
setting. We write these as the following
\begin{prop}
Let $M$ be (a domain of) a Riemannian manifold, $N$ a non-positive
curved manifold.
\\
1) (S-Y. Cheng \cite{ch})~Suppose that $u: M\to N$ is a harmonic
map with finite energy $E(u; M)$. Then on any compact subset $M'$
of $M$, one has the following estimate on energy density of $u$
\[
e(u)\le C(E(u), M', \dim M).
\]
2) (Schoen-Yau \cite{sy})~Let $u_1, u_2$ be two harmonic maps from
$M$ into $N$. Then the square of the distance function
$d_{N}^2(u_1, u_2)$ is subharmonic.
\end{prop}

We now show the following technical tool.

\begin{prop}
Let $K_1$ be an admissible metric with $K_1=h_1K_0$ under the
fixed flat basis, let $\mathcal{K}_0$ be the set of positive
self-adjoint operators $h$ on ${L_{\rho}}_{|\Delta^*_i}$ with
respect to $K_0$ satisfying that $h_{|r_i=1}={h_1}_{|r_i=1}$ and
$hK_0$ is admissible. Then the following functional
\begin{equation}\label{4.2}
G_i(h)=\int_{\Delta^*_i}\big(\mid(\delta'
h)h^{-1}\mid^2+\mid(\delta'' h)h^{-1}\mid^2\big)
\end{equation}
has a critical point $h$ on $\mathcal{K}_0$, satisfying
$$G_i(h)\le
G_i(h_1).
$$

Equivalently, there exists an admissible harmonic metric $K$
($=hK_0$) on ${L_\rho}_{|\Delta^*_i}$, which satisfies
$K_{|r_i=1}={K_1}_{|r_i=1}$ and
\[
\hat{E}_i(K)\le \hat{E}_i(K_1),
\]
where
$\hat{E}_i(K)=\int_{\Delta^*_i}\big(|\theta_K-\theta_{K_0}|^2+
|\overline\theta_K-\overline\theta_{K_0}|^2\big)$.

Furthermore, when considering $K$ and $K_0$ as equivariant maps
into $\mathcal{P}_n$, the maximum of the distance function
$d_{\mathcal{P}_n}(K(x), K_0(x))$ on $\Delta_i^*$ is attained on
$\{r_i=1\}$.
\end{prop}

\noindent {\it Remarks.} 1)~If $h$ is a function (e.g. the
representation is one-dimensional), just by setting $u=\log h$,
the problem is reduced to the usual Dirichlet problem for harmonic
functions on a disk. But the present situation is different, since
the covariant derivatives $\delta' h, \delta''h$ are defined using
the connections $\delta'+d'', \delta''+d'$ which are not defined
at the puncture $p_i$. In order to overcome this difficulty, we
use bounded exhausted domains with the fixed outside circle and
minimize the functional at each such domain; finally we show the
obtained sequence of minimizers converges to the required
minimizer of the functional. 2)~In the proof of Prop. 4, we only
use the harmonic property of $K_0$ on each $\Delta^*_i$.
\\
\\
{\it Proof of Proposition 4.} Take bounded exhausted domains
$\{\Delta_i^*\setminus\Delta_{i,{\frac 1k}}^*\}$ of $\Delta^*_i$,
where $\Delta^*_{i, {\frac{1}{k}}}$ represents the puncture disk
$\{t_i\in \mathbb{C}~|~0<|t_i|<1/k\}$. Minimizing the functional
on $\mathcal{K}_{0, k}=\{h\in \mathcal{K}_0~|~h|_{r_i={\frac
1k}}=h_1|_{r_i={\frac 1k}}\}${\footnote {More precisely, one
should put the functional on a closed convex subset of certain
Hilbert manifold.}}
\[
G_{i, {\frac 1k}}(h)=\int_{\Delta_i^*\setminus\Delta_{i, \frac
1k}^*}\big(\mid(\delta' h)h^{-1}\mid^2+\mid(\delta''
h)h^{-1}\mid^2\big),
\]
one gets a (unique) minimizer, denoted by $h_k$. This is
equivalent to minimizing the functional
\[
\hat{E}_{i, {\frac 1k}}(K)=\int_{\Delta^*_i\setminus\Delta_{i,
\frac 1k}^*}\big(|\theta_K-\theta_{K_0}|^2+
|\overline\theta_K-\overline\theta_{K_0}|^2\big)
\]
on the set $\{K=hK_0~|~h\in \mathcal{K}_{0, k}\}$, which is just a
boundary value problem for (equivariant) harmonic maps; since the
target space $\mathcal{P}_n$ is of non-positive sectional
curvature, this can always be solved uniquely by Hamilton (the
equivariant case by Schoen). We remark that the solution is a
(unique) minimizer of both $\hat{E}_{i, {\frac 1k}}(K)$ and
${E}_{i, {\frac 1k}}(K)=\int_{\Delta^*_i\setminus\Delta_{i, \frac
1k}^*}\big(|\theta_K|^2+ |\overline\theta_K|^2\big)$ under the
corresponding boundary condition, this is a direct consequence of
the Stokes formula and the harmonicity of $K_0$.
\\

We now show that the sequence $\{h_k\}$ (if necessary, going to a
subsequence) converges uniformly on any compact subset of
$\Delta^*_i$ to a critical point $h$ of $G_i$. The uniform
convergence of $h_k$ on compact subsets can be easily seen: Since
$h_k$ minimizes $G_{i, {\frac 1k}}$, they have uniform gradient
estimate in $k$ on any compact subset and hence $h_k$'s are
uniform bounded. The Arzela-Ascoli's Theorem then applies.
Actually one can further show $C^1$-convergence of $h_k$.
Furthermore, since $G_{i, {\frac 1k}}(h_k)\le G_{i, {\frac
1k}}(h_1)\le G_{i}(h_1)$, so one has
$$
G_i(h)\le G_{i}(h_1).
$$

Next, we need to show that that $h\in\mathcal{K}_0$, namely,
$K=hK_0$ and $K_0$ are mutually bounded; equivalently, this says
that when considering them as equivariant harmonic map into
$\mathcal{P}_n$, the distance function $d_{\mathcal{P}_n}(K, K_0)$
is bounded on $\Delta_i^*$. This can be obtained by using 2) of
Proposition 3 and Maximum Principle.

Applying Proposition 3, 2) to $d_{\mathcal{P}_n}^2(h_kK_0, K_0)$
on $\Delta_i^*\setminus\Delta^*_{i, {\frac 1k}}$, one has that for
all $k>1$, the distance functions $d_{\mathcal{P}_n}(h_kK_0, K_0)$
have $\max_{\Delta^*_i}d_{\mathcal{P}_n}(K_1, K_0)$ as an upper
bound. As $k$ goes to $\infty$, we obtain on $\Delta_i^*$
\[
d_{\mathcal{P}_n}(K,
K_0)\le\max_{\Delta^*_1}d_{\mathcal{P}_n}(K_1, K_0).
\]

We note that from the above argument, it is not very clear if $h$
(resp. $K$) is a minimizer of $G_i$ (resp. $\hat{E}_i$) on
$\mathcal{K}_0$; but from the following lemma, we will easily see
that this is actually the case.
\\

In order to prove the final assertion of Proposition 4, we first
state and prove the following
\begin{lem}
Let $K, K'$ be two harmonic metrics on $L_\rho|_{\Delta_i^*}$ with
the same boundary value $K_1|_{r_i=1}$ and satisfy that $K, K'$
are mutually bounded with $K_0$. Then $K\equiv K'$ on
$\Delta_i^*$.
\end{lem}
{\it Proof of Lemma 1.} The proof is obtained by again using
Porposition 3, 2) and the fact that on the half cylinder there
exists no nonconstant nonnegative bounded subharmonic function
which takes value zero on the boundary.~~~~~~~~~$\Box$
\\

\noindent{\it Continuation of Proof of Proposition 4.} Similar to
the argument in the beginning of the proof, we minimize the
functional
\[
{E}'_{i, {\frac 1k}}(K)=\int_{\Delta^*_i\setminus\Delta_{i, \frac
1k}^*}\big(|\theta_K|^2+ |\overline\theta_K|^2\big)
\]
on the set $\{K\in \mathcal{K}_0~|~K|_{r_i={\frac
1k}}=K_0|_{r_i={\frac 1k}}\}$. We obtain harmonic metrics
$h'_kK_0$ on $\Delta^*_i\setminus\Delta_{i, \frac 1k}^*$. Again by
Maximum Principle, $\max_{\Delta^*_i\setminus\Delta_{i, \frac
1k}^*}d_{\mathcal{P}_n}(h'_kK_0,
K_0)=\max_{r_i=1}d_{\mathcal{P}_n}(K_1, K_0)$. So, one can easily
show that these harmonic metrics converge to a harmonic metric
$K'=h'K_0$ on $\Delta_i^*$ and
$$
\max_{\Delta^*_i}d_{\mathcal{P}_n}(K',
K_0)=\max_{r_i=1}d_{\mathcal{P}_n}(K_1, K_0).
$$
Thus, by Lemma 1, we have $K\equiv K'$ and hence
$$
\max_{\Delta^*_i}d_{\mathcal{P}_n}(K,
K_0)=\max_{r_i=1}d_{\mathcal{P}_n}(K_1, K_0).
$$
This finishes the proof of Proposition 4.
~~~~~~~~~~~~~~~~~~~~~~~~~~~~$\Box$
\\

Let $K_1\in\mathcal{K}$ be an admissible metric on $L_\rho$.
Restricting $K_1$ to each $\Delta_i^*$ and using the solution of
Proposition 4 corresponding to ${K_1}_{|\Delta_i^*}$ to replace
${K_1}_{|\Delta_i^*}$, we obtain a new admissible metric on
$L_\rho$, denoted by $K$. Then, we have the following
\begin{prop}
\[
\hat{E}(K)\le\hat{E}(K_1).
\]
\end{prop}
{\it Proof.} By means of the definition of $K$ and Proposition 4,
we first have
\begin{eqnarray*}
&&\int_{X\setminus\cup_i\Delta^*_i}|\theta_K+\overline\theta_K|^2+
\sum_i\int_{\Delta^*_i}\big(|\theta_K-\theta_{K_0}|^2+
|\overline\theta_K-\overline\theta_{K_0}|^2\big)\\
&\le&\int_{X\setminus\cup_i\Delta^*_i}|\theta_{K_1}+\overline\theta_{K_1}|^2+\sum_i
\int_{\Delta^*_i}\big(|\theta_{K_1}-\theta_{K_0}|^2+
|\overline\theta_{K_1}-\overline\theta_{K_0}|^2\big).
\end{eqnarray*}
Note that $K\equiv K_1$ on ${X\setminus\cup_i\Delta^*_i}$. The
left-hand side the above inequality can be written as
\begin{eqnarray*}
\hat{E}(K)-\sum_i\int_{\Delta^*_i\setminus\Delta^*_{i/2}}|\theta_K+\overline\theta_K|^2+
\sum_i\int_{\Delta^*_i\setminus\Delta^*_{i/2}}\big(|\theta_K-\theta_{K_0}|^2+
|\overline\theta_K-\overline\theta_{K_0}|^2\big);
\end{eqnarray*}
similarly, the right-hand side is
\begin{eqnarray*}
\hat{E}(K_1)-\sum_i\int_{\Delta^*_i\setminus\Delta^*_{i/2}}|\theta_{K_1}+\overline\theta_{K_1}|^2+
\sum_i\int_{\Delta^*_i\setminus\Delta^*_{i/2}}\big(|\theta_{K_1}-\theta_{K_0}|^2+
|\overline\theta_{K_1}-\overline\theta_{K_0}|^2\big).
\end{eqnarray*}
So, in order to prove the Proposition, it suffices to prove the
following
\begin{eqnarray*}
&&-\sum_i\int_{\Delta^*_i\setminus\Delta^*_{i/2}}|\theta_K+\overline\theta_K|^2+
\sum_i\int_{\Delta^*_i\setminus\Delta^*_{i/2}}\big(|\theta_K-\theta_{K_0}|^2+
|\overline\theta_K-\overline\theta_{K_0}|^2\big)\\
&=&-\sum_i\int_{\Delta^*_i\setminus\Delta^*_{i/2}}|\theta_{K_1}+\overline\theta_{K_1}|^2+
\sum_i\int_{\Delta^*_i\setminus\Delta^*_{i/2}}\big(|\theta_{K_1}-\theta_{K_0}|^2+
|\overline\theta_{K_1}-\overline\theta_{K_0}|^2\big);
\end{eqnarray*}
this is equivalent to show
\[
\sum_i\int_{\Delta^*_i\setminus\Delta^*_{i/2}}<\theta_{K}+
\overline\theta_{K},
\theta_{K_0}+\overline\theta_{K_0}>=\sum_i\int_{\Delta^*_i\setminus\Delta^*_{i/2}}<\theta_{K_1}+
\overline\theta_{K_1}, \theta_{K_0}+\overline\theta_{K_0}>,
\]
which can be easily obtained by using the Stokes' formula together
with the facts that both $K$ and $K_1$ have the same boundary
value on each $\partial\Delta_i$ and that $K_0$ is harmonic on
each $\Delta^*_i$ and has vanishing normal derivative on
$\partial\Delta_{i/2}$.~~~$\Box$
\\

Now, we can turn to the proof of Main Theorem.
\\

\noindent {\it The proof of Main Theorem.} Let
$\{K_n\}_{n=1}^\infty$ be a minimizing sequence of $\hat{E}$ in
$\mathcal{K}$, i.e., $\lim_{n\to\infty}\hat{E}(K_n)=\mu$. Without
loss of generality, we can assume that {\it each metric $K_n$ is
harmonic on $X\setminus\cup_i\Delta^*_{i/2}$}; this can be done
briefly as follows: on $X\setminus\cup_i\Delta^*_{i/2}$, we
replace $K_n$ by a (unique) harmonic metric which is obtained by
solving the Dirichlet's boundary problem for equivariant harmonic
maps with the boundary value being
${K_n}_{|\cup_i\partial\Delta_{i/2}}$. Since this replacement does
not increase energy on $X\setminus\cup_i\Delta^*_{i/2}$, so the
new sequence is still a minimizing sequence.
\\

We now show that the minimizing sequence $\{K_n\}$ (if necessary,
going to a subsequence) is uniformly convergent on
$X\setminus\cup_i\Delta^*_{i}$ in the sense of $C^1$. The key is
to show the $C^0$-convergence; to this end, we use an idea due to
Donaldson \cite{do}. From now on, we consider each $K_n$'s as an
equivariant map from the universal covering of
$X\setminus\cup_i\Delta^*_{i/2}$ into $\mathcal{P}_m$. Since
$\{K_n\}$ is a minimizing sequence for $\hat{E}$, the usual energy
of $K_n$ on $X\setminus\cup_i\Delta^*_{i/2}$ is uniform bounded in
$n$. By means of Proposition 3, 1), the energy density $e(K_n)$ on
$X\setminus\cup_i\Delta^*_{i}$ has uniform bound in $n$, namely
\[
e(K_n)\le C, ~~{\text{on}}~X\setminus\cup_i\Delta^*_{i},
\]
where $C>0$ is independent of $n$.

Take a point $\tilde p$ in the universal covering of
$X\setminus\cup_i\Delta^*_{i}$, say, the projection of which lies
in $\cup_i\partial\Delta_{i}$, denoted by $p$. We would like to
show that the sequence $\{K_n(\tilde p)\}_{n=1}^\infty$ in
$\mathcal{P}_m$ (if necessary, going to a subsequence) converges.
If NOT, we may assume that $\{K_n(\tilde p)\}$ (if necessary,
going to a subsequence) converges to a point in a certain boundary
component, say $\Sigma$, of $\mathcal{P}_m$. By the
semi-simplicity of the representation $\rho: \pi_1(X)\to Gl(m,
\mathbb{C})$ (cf. \S3), there exists an element
$\gamma\in\pi_1(X)$ satisfying
$\rho(\gamma)\Sigma\cap\Sigma=\emptyset$. So, we have
\[
\lim_{n\to\infty}d_{\mathcal{P}_m}(K_n(\tilde p),
\rho(\gamma)K_n(\tilde
p))=\lim_{n\to\infty}d_{\mathcal{P}_m}(K_n(\tilde p),
K_n(\gamma\tilde p))=\infty.
\]
On the other hand, letting $c(t), t\in [0,1]$ be a differentiable
curve in the universal covering of $X\setminus\cup_i\Delta^*_{i}$
connecting the points $\tilde p$ and $\gamma(\tilde p)$, we then
have
\[
d_{\mathcal{P}_m}(K_n(\tilde p), K_n(\gamma\tilde p))\le
{\text{Length of the curve}}~ K_n(c(t)).
\]
By means of the uniform boundedness of energy density $e(K_n)$ in
$n$ on $X\setminus\cup_i\Delta^*_{i}$, we know that the length of
the curves $K_n(c(t))$ are actually uniformly bounded in $n$.
Thus, we derive a contradiction; namely, the sequence
$\{K_n(\tilde p)\}$ in $\mathcal{P}_m$ (if necessay, going to a
subsequence) converges.

Using the convergence of $\{K_n(\tilde p)\}$ and the fact that
$e(K_n)$ are uniformly bounded on $X\setminus\cup_i\Delta^*_{i}$,
we easily show that $K_n$ (if necessary, going to a subsequence)
is $C^1$-convergent on the compact subset
$X\setminus\cup_i\Delta^*_{i}$.
\\

Summing the above all up, we can assume that the minimizing
sequence $\{K_n\}$ in question satisfies that 1) on
$X\setminus\cup_i\Delta^*_{i/2}$, each $K_n$ is harmonic; 2) on
$X\setminus\cup_i\Delta^*_{i}$, $\{K_n\}$ uniformly converges in
the sense of $C^1$.
\\

Now, using Proposition 4 and 5, we construct a new minimizing
sequence from the above $\{K_n\}$. For each $K_n$, we restrict
$K_n$ to $\cup_i\Delta^*_i$ and consider this restriction as the
$K_1$ in Proposition 4 to get the corresponding $K$ in Proposition
4; we then use this $K$ to replace $K_n$ on $\cup_i\Delta^*_i$ to
get a new metric in $\mathcal{K}$, denoted by $K_n'$. Proposition
5 tells us $\hat{E}(K'_n)\le\hat{E}(K_n)$ so that the new sequence
$\{K_n'\}$ is still a minimizing sequence. Note that $K_n'\equiv
K_n$ on $X\setminus\cup_i\Delta^*_{i}$. Since $\{K_n\}$ (hence
$\{K_n'\}$) converges on $\cup_i\partial\Delta_i$, using the same
technique as in the proof of Proposition 4, we easily prove that
$\{K_n'\}$ converges on $\cup_i\Delta^*_{i}$; especially, the
limit is admissible. Thus, $\{K_n'\}$ converges on $X$ to an
admissible metric, denoted by $K$, and $K$ minimizes the modified
energy functional $\hat{E}$. By Proposition 2, we know that $K$ is
harmonic on $X$. ~~~~~~~~~~~~$\Box$

\vskip 1cm \noindent Department of Mathematics, Tongji University,
Shanghai, 200092, CHINA
\\
{\it E-mail}: {\verb"yhyang@mail.tongji.edu.cn"}

\end{document}